\documentclass{amsproc}
\newtheorem{theorem}{Theorem}[section]

\theoremstyle{definition}
\newtheorem{definition}[theorem]{Definition}
\newtheorem{example}[theorem]{Example}

\newtheorem{corollary}[theorem]{Corollary}
\theoremstyle{remark}
\newtheorem{remark}[theorem]{Remark}
\theoremstyle{theorem}
\newtheorem{proposition}[theorem]{Proposition}
\newtheorem{satz}{Theorem}
\numberwithin{equation}{section}

\newcommand{\C}{{\mathbb{C}}}
\newcommand{\D}{{\mathbb{D}}}
\newcommand{\R}{{\mathbb{R}}}
\newcommand{\N}{{\mathbb{N}}}

\renewcommand{\Re}{\mathop{\textrm{Re}}}
\def\BMOA{\mathop{BMOA}}

\begin{document}

\title[Critical points, the Gauss curvature equation \& Blaschke
  products]{Critical points, the Gauss curvature equation and\\ Blaschke
  products}

\author{Daniela Kraus}
\address{Department of Mathematics, University of W\"urzburg, Emil--Fischer
  Stra{\ss}e 40, 97074 W\"urzburg, Germany}
\email{dakraus@mathematik.uni-wuerzburg.de}
\thanks{Supported by the Deutsche Forschungsgemeinschaft (Grants: Ro 3462/3--1 and Ro 3462/3--2)}

\author{Oliver Roth}
\address{Department of Mathematics, University of W\"urzburg, Emil--Fischer
  Stra{\ss}e 40, 97074 W\"urzburg, Germany}
\email{roth@mathematik.uni-wuerzburg.de}

\date{January 1, 1994 and, in revised form, June 22, 1994.}

\begin{abstract}
In this survey paper we discuss the problem of characterizing the critical
sets of bounded analytic functions in the unit disk of the complex plane. 
This problem is closely
related to the Berger--Nirenberg problem in differential geometry as well as
to the
problem of describing the zero sets of functions in Bergman spaces.
It turns out that for any non--constant bounded analytic function in the unit
disk there is always a (essentially) unique ``maximal'' Blaschke product with the same
critical points. These  maximal Blaschke products  have remarkable 
properties simliar to those of Bergman space inner functions and they provide a natural
generalization of the class of finite Blaschke products.
\end{abstract}

\maketitle

\section{Critical points of bounded analytic functions}

A sequence of points $(z_j)$ in a subdomain $\Omega$ of the complex plane $\C$ is called the
zero set of an analytic function $f : \Omega \to \C$, if $f$ vanishes precisely on this
set. This means that if the point $\xi\in \Omega$ occurs $m$ times in the sequence, then $f$
has a zero at $\xi$ of precise order $m$, and $f(z)\not=0$ for $z \in \Omega
\backslash (z_j)$. 
The following classical theorem 
due to  Jensen \cite{Jen1899},
Blaschke \cite{Bla1915} and F.~and R.~Nevanlinna \cite{Nevs1922}  characterizes completely the zero
sets of bounded analytic functions defined in the open unit disk $\D:=\{z \in
\C \, : \, |z|<1\}$.

 \begin{satz}\label{thm:classical_zero}
 Let $ (z_j)$ be a sequence in  $\D$. Then the following statements are
equivalent. 
\begin{itemize}
 \item[(a)] There is an analytic self--map of $\D$ with zero set
 $ (z_j)$.
\item[(b)]  There is a Blaschke product with zero set  $
  (z_j)$.
\item[(c)] The sequence $(z_j)$ fulfills the Blaschke condition, i.\!\;e.~$
\sum\limits_{j=1}^{\infty} \big(1-|z_j|\big) < + \infty\, .
$
\end{itemize}
\end{satz}

We call a sequence of points $(z_j)$ in a domain $\Omega \subseteq \C$ the
critical set of an analytic function $f : \Omega \to \C$, if $(z_j)$ is the zero
set of the first derivative $f'$ of the function $f$.
 There is an extensive literature on critical sets.
In particular, there are many interesting results on the relation
between the zeros and the critical points of analytic and harmonic functions.
 A classical reference for all this is the book of Walsh \cite{Wal1950}.

The first aim of this survey paper is to point out the following analogue of
Theorem~\ref{thm:classical_zero} for the critical sets of bounded analytic
functions instead of their zeros sets.

\begin{theorem}\label{thm:main}
Let $ (z_j)$ be a sequence in $\D$. Then the following statements are
equivalent.
\begin{itemize}
 \item[(a)] There is an analytic  self--map of $\D$ with critical set
 $ (z_j)$.
\item[(b)]  There is a Blaschke product with critical set  $
  (z_j)$.
\item[(c)] There is  a function in the weighted Bergman space $$ \mathcal{
    A}^2_1=\left\{
                  g:\D \to \C \text{ analytic} \,:  \iint \limits_{\D} (1-|z|^2)\, |g(z)|^2\,
d\sigma_z <+ \infty \right\} $$ with zero set  $ (z_j)$.
\end{itemize}
\end{theorem}

Here, and throughout, $\sigma_z$ denotes two--dimensional Lebesgue measure
with respect to $z$.

For the special case of  {\it finite} sequences, a result related to Theorem
\ref{thm:main}  can be found in  work of
Heins \cite[\S29]{Hei62}, 
Wang \& Peng \cite{WP79}, Zakeri \cite{Z96} and  Stephenson
\cite{Ste2005}. They proved that for 
every finite sequence $(z_j)$ in $\D$ there is always a
 {\it finite} Blaschke product whose critical set coincides with $(z_j)$, see
 also Remark \ref{rem:finite} below. 
A recent  generalization of this result to {\it infinite} sequences 
is discussed in \cite{KR}. There it is shown that every Blaschke sequence
$(z_j)$ is the critical set  of an infinite Blaschke product. However,
the converse to this, known as the Bloch--Nevanlinna conjecture
\cite{Dur1969}, is false. According to a result of Frostman, there do exist
 Blaschke products whose critical sets
fail to satisfy the Blaschke condition, see \cite[Theorem 3.6]{Col1985}. 
Thus the critical sets of bounded analytic functions are not
 just the Blaschke sequences and the situation for critical sets 
is much more complicated than for  zero sets.

  Theorem \ref{thm:main} identifies the critical
 sets of bounded analytic functions as the zeros sets of functions in the
 Bergman space~$\mathcal{ A}^2_1$.
The simple geometric characterization of the zero sets of bounded analytic
functions via the Blaschke condition (c) in Theorem \ref{thm:classical_zero} has
not found an explicit counterpart for 
critical sets yet.  However, condition (c) of Theorem \ref{thm:main} might be seen as
an implicit substitute.
The zero sets of (weighted) Bergman space functions have intensively been studied in the 1970's and
1990's by Horowitz
\cite{Hor1974,Hor1977}, Korenblum \cite{Kor1975}, Seip \cite{Seip1993,
  Seip1995}, Luecking \cite{Lue1996} and many others. As a result quite sharp necessary conditions and  sufficient
conditions for the zero sets of  Bergman
space functions  are available. 
In view of Theorem \ref{thm:main} all results about zero sets of
Bergman space functions carry now
over to the critical sets of bounded analytic functions and vice versa.
Unfortunately, a {\it geometric}
characterization of the zero sets of  (weighted) Bergman space functions is
still unknown, ``and it is well known
that this problem is  very difficult'', cf.~\cite[p.~133]{HKZ}.

The implications ``(b) $\Longrightarrow$ (a)'' and ``(a) $\Longrightarrow$ (c)''
of Theorem \ref{thm:main} are easy to prove. In fact, the statement ``(a)
$\Longrightarrow$ (c)''  follows directly from the Littlewood--Paley identity (see
\cite[p.~178]{Sha1993}), which says that for any holomorphic function $f : \D \to \D$ the derivative
$f'$ belongs to $\mathcal{A}^2_1$. Hence the critical set of any non--constant
bounded analytic function is trivially the zero set of a function in
$\mathcal{A}^2_1$. It is however not true that any $\mathcal{A}^2_1$ function
is the derivative of a bounded analytic function, so ``(c) $\Longrightarrow$
(a)'' is more subtle.

In the following sections, we describe some of the ingredients of the proofs of
the implications ``(a) $\Longrightarrow$ (b)'' and ``(c) $\Longrightarrow$
(a)'' in
Theorem \ref{thm:main} as well as some further results connected to it. 
In particular, we explain the close relation of Theorem~\ref{thm:main}
to conformal differential geometry and the solvability of the Gauss curvature equation.
This paper is expository, so there are  essentially no proofs. 
For the proofs we refer to \cite{Kra2010, Kra2011a, KRR06,KR, KR2011}.
Background material on Hardy spaces and Bergman spaces can be found e.g.~in
the excellent books \cite{Dur2000,DS,Gar2007,HKZ,Koo1998,Mas2009}.

Finally, we draw attention to the recent paper \cite{EKS2012} of P.~Ebenfelt,
D.~Khavinson and H.S.~Shapiro and the references therein, where the difficult
problem of constructing  a finite Blaschke product by its \textit{critical
  values} is discussed.

\section{Conformal metrics and maximal Blaschke products} \label{sec:2}

An essential characteristic of the proof of  
the implication ``(a) $\Longrightarrow$ (b)'' in Theorem \ref{thm:main} is the extensive use of
negatively curved conformal pseudometrics. We give a short account of some of
their properties and refer to \cite{BM2007, Hei62, KL2007, Krantz, KR2008, Smi1986} for more details.
In the following, $G$ and $D$
always denote domains in the complex plane $\C$.

\subsection{Conformal metrics  and developing maps}

We recall that a nonnegative upper
semicontinuous function $\lambda$ on $G$, $\lambda: G \to [0, + \infty)$, $\lambda \not\equiv 0$, is
called conformal density on $G$.  The corresponding
 quantity \label{def:leng}
$\lambda(z) \, |dz|$   is called conformal pseudometric on $G$.
If $\lambda(z)>0$ for all $z \in G$, we say $\lambda(z) \, |dz|$ is a
conformal metric on $G$. We call a conformal
pseudometric $\lambda(z) \, |dz|$  regular on $G$, if $\lambda$ is of class
$C^2$ in $\{z \in G \, : \, \lambda(z)>0\}$.

\begin{example}[The hyperbolic metric]
The ubiquitous example of a conformal metric is the Poincar\'{e} metric or hyperbolic
metric $$\lambda_{\D}(z)\, |dz|:=\frac{|dz|}{1-|z|^2}$$ for
 the unit disk $\D$. Clearly, $\lambda_{\D}(z)\, |dz|$  is a regular conformal
 metric on $\D$.
\end{example}

The (Gauss) curvature $\kappa_{\lambda}$ of a regular conformal pseudometric
$\lambda(z) \, |dz|$ on $G$
is defined by
\begin{equation*}
\kappa_{\lambda}(z):=-\frac{\Delta (\log \lambda)(z)}{\lambda(z)^2}
\end{equation*}
 for all points $z \in G$ where $\lambda(z) > 0$.
Note that if $\lambda(z)\, |dz|$ is a regular conformal metric
with curvature $\kappa_{\lambda}=\kappa$ on $G$, then the function $u:= \log \lambda$ 
satisfies the partial differential equation (Gauss curvature equation)
\begin{equation}\label{eq:pde1} 
\Delta u=-\kappa(z) \, e^{2 u}
\end{equation}
on $G$. If, conversely, a realvalued $C^2$--function $u$ is a solution of the 
Gauss curvature equation (\ref{eq:pde1}) on $G$,
then $\lambda(z):=e^{u(z)}$ induces a regular conformal metric $e^{u(z)} \
|dz|$ on $G$ with
curvature $\kappa$.

\begin{example}
The hyperbolic metric $\lambda_{\D}(z) \, |dz|$ on $\D$ has constant
negative curvature\footnote{\textbf{Warning:} A rival school of thought calls $\frac{2 \,
    |dz|}{1-|z|^2}$ the hyperbolic metric of $\D$. This metric has constant
  curvature $-1$.}  $-4$.
\end{example}

Using analytic maps, conformal metrics can be transferred from one domain to
another as follows. Let $\lambda(w) \, |dw|$ be a conformal pseudometric
on a domain $D$  and let $w=f(z)$ be  a non--constant analytic map from another
domain $G$ to $D$. Then
the conformal pseudometric
\begin{equation*} 
 (f^*\lambda)(z) \, |dz|:=\lambda(f(z)) \, |f'(z)| \, |dz|\, , 
\end{equation*}
defined on $G$, is called the pullback of $\lambda(w) \, |dw|$ under the map
$f$. 
Now, Gauss curvature is important since it is a conformal invariant in
the following sense.

\begin{satz}[Theorema Egregium]
For every analytic map $w=f(z)$ and every regular conformal pseudometric $\lambda(w)\,|dw|$ the relation
\begin{equation*} 
\kappa_{f^*\lambda}(z)=\kappa_{\lambda}(f(z))
\end{equation*}
is satisfied provided $\lambda(f(z)) \,  |f'(z)|>0$. 
\end{satz}

Hence pullbacks of  conformal pseudometrics can be used for constructing conformal
pseudometrics with various prescribed properties while controlling their
curvature. For instance,  if $\lambda(w) \, |dw|$ is a conformal metric on $D$
(without zeros) and $w=f(z)$ is a non--constant analytic map from $G$ to $D$,
then $(f^*\lambda)(z) \, |dz|$ is a pseudometric on $G$ with zeros exactly at the
critical points of $f$. In order to take  the multiplicity of the
zeros into account, it is convenient to introduce the 
 following  formal definition. 

 \begin{definition}[Zero set of a pseudometric]\label{def:zeroset}
Let $\lambda(z)\, |dz|$ be a conformal pseudometric on $G$. We say
$\lambda(z)\, |dz|$ has a zero of order $m_0>0$ at $z_0 \in G$ if
\begin{equation*}
\lim_{z \to z_0} \frac{\lambda(z)}{|z-z_0|^{m_0}} \quad\text{ exists and }
\not=0 \, . 
\end{equation*}
We will call a  sequence $\mathcal{ C}=(\xi_j) \subset G$
\begin{equation*}
(\xi_j):=(\underbrace{z_1, \ldots, z_1}_{m_1  -\text{times}},\underbrace{z_2,
  \ldots, z_2}_{m_2-\text{times}} , \ldots  )\,, \, \,  z_k \not=z_n \text{ if
} k\not=n, 
\end{equation*}
  the zero set of a conformal
pseudometric $\lambda(z) \, |dz|$, if $\lambda(z)>0$ for $z \in
G\backslash \mathcal{ C}$ and if $\lambda(z)\, |dz|$ has a zero of order $m_k \in \N$
at $z_k$ for all $k$.
\end{definition}

Now, every conformal pseudometric of constant curvature $-4$
can  locally be
represented by a holomorphic function. This is the content of the following result.

\begin{satz}[Liouville's Theorem]
\label{thm:liouville}

Let $\mathcal{ C}$ be a sequence of points in a simply connected domain $G$
and let
$\lambda(z) \, |dz|$ be a regular conformal pseudometric on $G$ with constant curvature
$-4$ on $G$ and zero set $\mathcal{ C}$. Then $\lambda(z) \, |dz|$ is the pullback of the
hyperbolic metric $\lambda_{\D}(z)\, |dz|$ under some analytic map $f:G \to
\D$, i.\!\;e.
\begin{equation}\label{eq:liouville}
 \lambda(z) =\frac{|f'(z)|}{1-|f(z)|^2}\, , \quad z \in G. 
\end{equation}
If $g:G  \to \mathbb{D}$ is another analytic function, then
$\lambda(z)=(g^*\lambda_{\D})(z)$ for $z\in G$
 if and only if $g=T\circ f$ for some conformal automorphism $T$ of $\D$. 
\end{satz}


A holomorphic function $f$ with  property (\ref{eq:liouville})
 will be called a {\bf
  developing map} for $\lambda(z) \, |dz|$. 
Note that the critical set of each developing map coincides with the zero set
of the corresponding conformal pseudometric.  

\begin{example}
The developing maps of the hyperbolic metric $\lambda_{\D}(z) \, |dz|$
are precisely  the conformal automorphisms of $\D$, i.e., the finite Blaschke
products of degree $1$.
\end{example}

For later applications we wish to mention the following variant of Liouville's Theorem.

\begin{theorem} \label{rem:liouville}
Let $G$ be a simply connected domain and let $h:G \to \C$ be an analytic
map. If $\lambda(z)\, |dz|$ is a regular conformal metric with curvature $-4
\, |h(z)|^2$, then there exists a holomorphic function $f:G \to \D$ such
that
\begin{equation*}
 \lambda(z) =\frac{1}{|h(z)|}\frac{|f'(z)|}{1-|f(z)|^2}\, , \quad z
\in G. 
\end{equation*}
Moreover, $f$ is uniquely determined up to postcomposition with a unit disk automorphism.
\end{theorem}

 Liouville \cite{Lio1853} stated
Theorem \ref{thm:liouville}   for the special case that $\lambda(z)
\, |dz|$ is a regular conformal metric.
We therefore  refer to
 Theorem \ref{thm:liouville} as well as to Theorem \ref{rem:liouville} as
 Liouville's theorem. 
 Theorem \ref{thm:liouville} and in particular the special case that $\lambda(z) \, |dz|$ is
 a conformal metric has a number of different proofs, see for instance 
\cite{Bie16, CW94, CW95, Min, Nit57, Yam1988}. Theorem
\ref{rem:liouville} is discussed in \cite{KR}.

Liouville's theorem plays an important r$\hat{\text{o}}$le in this paper. It
provides a bridge bet\-ween the world of bounded analytic functions and the
world of conformal pseudometrics with constant negative curvature.
The critical points on the one side correspond to the zeros on the other side.
Unfortunately, this bridge only works for simply connected domains, see Remark \ref{rem:final}.

\subsection{Maximal conformal pseudometrics  and maximal Blaschke products}

Apart from having constant negative curvature the hyperbolic metric 
on $\D$ has another important property: it
 is \textit{maximal} among all regular conformal pseudometrics on
$\D$ with curvature bounded above by $-4$. This is the content of the
following result.

\begin{satz}[Fundamental Theorem]\label{thm:ahlfors}
Let $\lambda(z)\,|dz|$ be a regular conformal pseudometric on
$\D$ with curvature bounded above by $-4$. Then $\lambda(z) \le
\lambda_{\D}(z)$ for every $z \in \D$.
\end{satz}

Theorem \ref{thm:ahlfors} is due to Ahlfors \cite{Ahl1938} and it is usually  
called Ahlfors' lemma. However, in view of its relevance Beardon and Minda
proposed to call Ahlfors' lemma the  {\bf fundamental theorem}. We
will follow  their suggestion in this paper. As a result of the fundamental
theorem we  have
\begin{eqnarray*}
 \lambda_{\D}(z)=\max \{ \lambda(z) \, : \, \lambda(z)\,|dz| \text{ is a regular conformal pseudometric on
$\D$} & &  \\ & & \hspace*{-4cm}  \text{ with curvature }  \le -4  \}
\end{eqnarray*}
for any $z \in \D$.

\begin{remark}[Developing maps and universal coverings] \label{rem:hyp_all}
Let $G$ be a hyperbolic subdomain of the complex plane $\C$, i.e., the
complement $\C\backslash G$ consists of more than one point.
In analogy with the Poincar\'e metric for the unit disk, a regular conformal metric
$\lambda_{G}(z)\, |dz|$ of constant curvature $-4$ on  $G$ is said to be the
hyperbolic metric for $G$, if it is the maximal regular conformal pseudometric with
curvature $\le -4$ on $G$, i.\!\;e.~$\lambda(z) \le \lambda_G(z)$, $z \in G$, for all
regular conformal pseudometrics $\lambda(z)\, |dz|$ on $G$ with curvature
bounded above by $-4$. Then
$\lambda_G(z)\, |dz|$ and $\lambda_{\D}(z)\, |dz|$ are
connected via the universal coverings $\pi:\D \to G$ of $G$ by the formula
$(\pi^*\lambda_{G})(z)\, |dz|=\lambda_{\D}(z)\, |dz|$,
see for example \cite[\S 9]{Hei62} and \cite{KR2008, Min1979}.
Hence,  every branch of the inverse of a universal covering map $\pi : \D \to G$
is \textit{locally} the developing map of the hyperbolic metric $\lambda_G(z)
\, |dz|$ of $G$.
In particular, if $G$ is a hyperbolic simply connected domain, then the
developing maps for $\lambda_{G}(z) \, |dz|$ are precisely the conformal
mappings from $G$ onto $\D$.
\end{remark}

We now consider prescribed zeros.

\begin{satz} \label{thm:perron}
Let $\mathcal{C}=(\xi_j)$ be a sequence of points in $G$  and 
\begin{eqnarray*}
\Phi_{\mathcal{C}}:=\left\{ \lambda  \, : \, \lambda(z) \, |dz| \text{ is a regular conformal
    pseudometric in } G \right. & & \\ & & \hspace*{-5.8cm} \left. \text{ with
    curvature } \le -4 \text{ and zero set } \mathcal{C}^*
  \supseteq \mathcal{C}\right\}\, .
\end{eqnarray*}
If $\Phi_{\mathcal{C}}\not=\emptyset$,  then
$$ \lambda_{max}(z):=\sup\{ \lambda(z) \, : \, \lambda \in \Phi_{\mathcal{C}}\} \, , \qquad
z \in G \, , $$
induces the unique maximal regular conformal pseudometric
$\lambda_{max}(z) \, |dz|$ on $G$ with constant
curvature $-4$ and zero set $\mathcal{C}$.
\end{satz}

Thus, if $\Phi_{\mathcal{C}}\not=\emptyset$, i.e.,
if there exists at least one regular conformal pseudometric $\lambda(z) \, |dz|$
on $G$ with curvature $\le -4$ whose zero set contains the given sequence $\mathcal{C}$,
then there exists a (maximal) regular conformal pseudometric $\lambda(z) \, |dz|$
on $G$ with \textit{constant} curvature $-4$ whose zero set is
\textit{exactly} the sequence $\mathcal{C}$. 
In particular, Theorem \ref{thm:perron} can be applied,  if
there exists a non--constant holomorphic function $f : G \to \D$ with critical
set $\mathcal{C}^* \supseteq \mathcal{C}$ since then the pseudometric $(f^*\lambda_{\D})(z) \, |dz|$
belongs to $\Phi_{\mathcal{C}}$.
The proof of Theorem \ref{thm:perron}
relies on a modification of Perron's method and can be found in \cite[\S 12 \&
\S 13]{Hei62}
and \cite{Kra2011a}.

\begin{example}
If $G$ is a hyperbolic domain and $\mathcal{C}=\emptyset$, then the maximal
regular conformal pseudometric $\lambda_{max}(z) \, |dz|$ on $G$ with constant
curvature $-4$ and zero set $\mathcal{C}$ is exactly the hyperbolic metric $\lambda_G(z)
\, |dz|$ for $G$.
\end{example}

Thus maximal pseudometrics are generalizations of the hyperbolic metric and
their developing maps are therefore of special interest.

\begin{definition}[Maximal functions]
Let $\mathcal{C}$ be a sequence of points in $\D$  such that there exists a  maximal regular conformal pseudometric
$\lambda_{max}(z) \, |dz|$ on $\D$ with constant
curvature $-4$ and zero set $\mathcal{C}$. Then every developing map for
$\lambda_{max}(z) \, |dz|$ is called maximal function for $\mathcal{C}$.
\end{definition}

Some remarks are in order. First, 
in view of Theorem \ref{thm:liouville} every maximal function is uniquely
determined by its critical set $\mathcal{C}$ up to postcomposition with a unit disk automorphism and, conversely,
the postcomposition of any maximal function with a unit disk
automorphism is again a maximal function. Second, if  $\mathcal{C}=\emptyset$,
then the maximal functions for $\mathcal{C}$ are precisely the unit disk
automorphisms, i.e., the finite Blaschke products of degree~$1$.

Now, we have the following result, see \cite{Kra2011a}.

\begin{theorem} \label{thm:0}
Every maximal function  is a Blaschke product.
\end{theorem}

It is to emphasize that
 since the postcomposition of any maximal function with a unit disk
automorphism is again a maximal function, every maximal function is an
\textit{indestructible} Blaschke product.

We note that Theorem \ref{thm:perron} combined with Theorem \ref{thm:0}
immediately gives the
implication ``(a) $\Longrightarrow$ (b)'' in Theorem \ref{thm:main}.
In fact, if $f$ is a non--constant analytic self--map of $\D$ with critical
set $\mathcal{C}$, then $\lambda(z) \, |dz|:=(f^*\lambda_{\D})(z) \, |dz|$
is a regular conformal pseudometric on $\D$ with zero set $\mathcal{C}$. Thus
Theorem \ref{thm:perron} guarantees the existence of a maximal conformal
pseudometric on $\D$ with zero set $\mathcal{C}$. Now, Theorem \ref{thm:0}
says that the corresponding maximal function for $\mathcal{C}$ is a Blaschke product.

In the special case that the maximal function has finitely many critical
points, the statement of Theorem \ref{thm:0} follows from the next result which is due to Heins 
\cite[\S 29]{Hei62}.


\begin{satz} \label{satz:heins}
Let $\mathcal{C}=(z_1, \ldots, z_n)$ be a finite sequence in $\D$
and $f : \D \to \D$ analytic.
 Then the following statements are equivalent.
\begin{itemize}
\item[(a)] $f$ is a maximal function for $\mathcal{C}$.
\item[(b)] $f$ is Blaschke product of degree $n+1$ with critical set $\mathcal{C}$.
\end{itemize}
\end{satz}

We shall give a quick proof of Theorem \ref{satz:heins} in Remark
\ref{satz:heinsproof} below.

\begin{remark}[Constructing finite Blaschke products with prescribed critical
  points] \label{rem:finite}
In his proof of  Theorem \ref{satz:heins}, Heins showed that for any finite sequence
$\mathcal{C}=(z_1, \ldots, z_n)$ in $\D$  there is always a finite Blaschke
product $B$ of
degree $n+1$ with critical set $\mathcal{C}$. The essential step is
nonconstructive and consists in showing that the set of critical points of all
finite Blaschke products of degree $n+1$, which is clearly closed, is also open in the
poly disk $\D^n$ by applying Brouwer's fixed point theorem.
The same result was later obtained by  Wang \& Peng \cite{WP79} and Zakeri
 \cite{Z96} by using similar arguments. A completely different approach 
 via Circle Packing is due to Stephenson, see Lemma 13.7 and Theorem 21.1 in \cite{Ste2005}. Stephenson builds 
 discrete finite Blaschke products with prescribed branch set and shows that
 under refinement  these
 discrete Blaschke products converge locally uniformly in $\D$ to the desired  classical
 Blaschke product.

We are not aware of any
 \textit{efficient} constructive method for computing a (nondiscrete) finite Blaschke
 product from its critical points.
\end{remark}

Maximal functions form a particular class of Blaschke products. It is
therefore convenient to make the following definition.

\begin{definition}[Maximal Blaschke products]
A non--constant Blaschke product  is called a maximal Blaschke product, if it is a
maximal function for its critical set.
\end{definition}

As was mentioned earlier,
 every maximal Blaschke product is indestructible and every finite
Blaschke product is a maximal Blaschke product. Moreover,  if $\mathcal{C}$ is the
critical set of any non--constant analytic function $f : \D \to \D$, then there
is maximal Blaschke product with critical set $\mathcal{C}$. A maximal Blaschke
product is uniquely determined by its critical set $\mathcal{C}$ up to postcomposition with a unit disk automorphism.

Geometrically,
 the finite maximal Blaschke products  are just the finite branched
coverings of $\D$. One is therefore inclined to consider maximal Blaschke
products  for
infinite branch sets as ``infinite branched coverings'':

\renewcommand{\arraystretch}{1.5}
\begin{center}
\begin{tabular}{|l|p{4.7cm}|p{4.8cm}|}
\hline
 {\bf critical set} & { \bf maximal Blaschke product} & {\bf mapping properties }\\[2mm] \hline
$\mathcal{ C}=\emptyset$ & automorphism of $\D$& unbranched covering of $\D$\!\;; \\[-1.5mm]
&&conformal self--map of $\D$\\[0.5mm] \hline
$\mathcal{ C}$ finite &finite Blaschke product &finite branched covering of $\D$\\[0.5mm] \hline
$\mathcal{ C}$ infinite & indestructible infinite maximal Blaschke product& ``\!\;infinite
branched covering of $\D$\!\;''\\ \hline
\end{tabular}
\end{center}

The class of maximal conformal pseudometrics and their corresponding maximal
functions have already been studied by Heins in \cite[\S25 \& \S
26]{Hei62}. Heins obtained
  some necessary conditions as well as sufficient conditions   
for maximal functions (see Theorem \ref{thm:omitted_value} and Theorem
\ref{thm:heins_island} below), but he did not prove that maximal functions are always
Blaschke products. He also posed the 
problem of characterizing maximal functions, cf.~\cite[\S 26 \& \S 29]{Hei62}.


\section{Some properties of maximal Blaschke products}

It turns out that maximal Blaschke products do have remarkable 
properties and provide in some sense a fairly natural generalization of the class of finite
Blaschke products. In this section we 
take a closer look at some of the properties of maximal Blaschke
products. In the following, we denote by  $H^{\infty}$ the set of bounded analytic functions $f : \D \to
\C$ and set $||f||_{\infty}=\sup_{z \in \D} |f(z)|$. This makes
$(H^{\infty},||\cdot||_{\infty})$ a Banach space.

\subsection{Schwarz' lemmas}

\begin{theorem}[Maximal Blaschke products as extremal functions] \label{thm:2}
Suppose $\mathcal{ C}$ is a sequence of
points in $\D$ such that
\begin{equation*}
 \mathcal{ F}_{\mathcal{ C}}:=\big \{ f \in H^{\infty} \, : \,  
f'(z)=0 \text{ for } z \in \mathcal{ C} \big\}
\end{equation*}
contains at least some non--constant function and let  $$ m:=\min\{\,  n \in \N \, : \, f^{(n)}(0)\not=0 \text{ for some } f \in \mathcal{
  F}_\mathcal{ C} \, \}\ge 1\, .$$ Then  the unique extremal function to the extremal problem
\begin{equation} \label{eq:ex}
 \max \big\{ \Re f^{(m)}(0) \, : f \in \mathcal{ F}_\mathcal{ C}, \,
 {||f||}_{\infty}\le 1\big \}
\end{equation}
is  the maximal Blaschke product  $F$ for $\mathcal{ C}$
normalized by $F(0)=0$ and $F^{(m)}(0)>0$.  
\end{theorem}

We refer to \cite{KR2011} for the proof of Theorem \ref{thm:2}. To put
Theorem \ref{thm:2} in perspective, note that if $m=1$, then the extremal
problem (\ref{eq:ex}) is exactly the problem of maximizing the derivative at a
point, i.e., exactly the character of Schwarz' lemma.

\begin{remark}[The Nehari--Schwarz Lemma]

If $\mathcal{C}=\emptyset$, then
  $\mathcal{F}_{\mathcal{C}}=H^{\infty}$, i.e., the set of all bounded analytic
  functions in $\D$. In this case, Theoren \ref{thm:2} is of course just the statement of Schwarz' lemma.
If $\mathcal{C}$ is a finite sequence, then  Theorem \ref{thm:2} is exactly
  Nehari's
1947 generalization of Schwarz' lemma (Nehari \cite{Neh1946}, in particular the
Corollary to Theorem 1). Hence Theorem \ref{thm:2} can be considered as  an extension of the
Nehari--Schwarz Lemma.
\end{remark}

\begin{remark}[The Riemann mapping theorem and the Ahlfors map]

We consider a domain  $\Omega \subseteq \C$  containing $0$. Let $\mathcal{C}=(z_j)$
  be the  critical set of a non--constant function $f$ in $H^{\infty}(\Omega)$,
  where $H^{\infty}(\Omega)$ denotes the set of all functions analytic and
  bounded in $\Omega$. We let $N$ denote the number 
of times that $0$ appears in the sequence~$\mathcal{  C}$ and set
$$ \mathcal{F}_{\mathcal{C}}(\Omega):=\left\{f \in H^{\infty}(\Omega) \, : \, f'(z)=0 \text{ for
  any } z \in \mathcal{C}\right\} \, .$$
Then, by a normal family argument, there is always at least one extremal
function for the extremal problem
\begin{equation} \tag{$*$}
 \max \big\{ \Re f^{(N+1)}(0) \, : f \in \mathcal{ F}_\mathcal{ C}(\Omega), \, {||f||}_{\infty} \le 1 \big \} \, .
\end{equation}
In the following three cases there is a unique extremal function to the
extremal problem ($*$).
\begin{itemize}
\item[(i)] \textit{$\Omega\not=\C$ is simply connected and
  $\mathcal{C}=\emptyset$ (conformal maps):}\\
 In this case, $N=0$ and the
  extremal problem ($*$) has exactly one extremal function,
  the normalized Riemann map
  $\Psi$ for $\Omega$, that is, the unique conformal map $\Psi$ from $\Omega$
  onto $\D$ normalized such that   $\Psi(0)=0$ and $\Psi'(0)>0$.
 \item[(ii)] \textit{$\Omega \not= \C$ is simply connected and
     $\mathcal{C}\not=\emptyset$ (prescribed critical points):}\\
Let $\Psi$ be the normalized Riemann map for $\Omega$.
Then $\Psi(\mathcal{C})$ is the critical set of a non--constant function in
$H^{\infty}=H^{\infty}(\D)$. If 
$B_{\Psi(\mathcal{C})}$ is the extremal function in
$\mathcal{F}_{\Psi(\mathcal{C})}$ according to Theorem \ref{thm:2}, then
$B_{\Psi(\mathcal{C})} \circ \Psi$ is the unique extremal function for ($*$).
\item[(iii)] \textit{$\Omega \not= \C$ is not simply connected  and
  $\mathcal{C}=\emptyset$ (Ahlfors' maps):}\\
If $\Omega$ has connectivity $n \ge 2$, none of whose boundary components
reduces to a point, then the
  extremal problem ($*$) has exactly one solution, namely the Ahlfors map
  $\Psi : \Omega \to \D$. It is a $n:1$ map  from $\Omega$ onto $\D$ such  that 
  $\Psi(0)=0$ and $\Psi'(0)>0$, see Ahlfors \cite{Ahl} and Grunsky \cite{Gru1978}.
\end{itemize}
\end{remark}

The   Schwarz lemma (i.e., the case $\mathcal{C}=\emptyset$ of Theorem \ref{thm:2})
can be stated in  an invariant form,  the Schwarz--Pick lemma which says that
\begin{equation} \label{eq:sp}
 \frac{|f'(z)|}{1-|f(z)|^2} \le \frac{1}{1-|z|^2} \, , \qquad z \in \D \, , 
\end{equation}
for any analytic map $f : \D \to \D$, with equality for some point $z \in \D$
if and only if $f$ is a conformal disk
automorphism. Hence maximal Blaschke products without critical points
serve as extremal functions. In a similar way, the more general statement of Theorem
\ref{thm:2} admits an invariant formulation as follows (see \cite{Kra2011a}).

\begin{theorem}[Sharpened Schwarz--Pick inequality]\label{cor:3}
Let $f: \D \to \D$ be a non--constant analytic function with critical set $\mathcal{
  C}$ and let $\mathcal{ C}^*$ be a subsequence of $\mathcal{ C}$. 
Then there exists a maximal Blaschke
product $F$ with critical set $\mathcal{ C}^*$
such that
\begin{equation*}
\frac{|f'(z)|}{1-|f(z)|^2} \le \frac{|F'(z)|}{1-|F(z)|^2}\, ,\quad  z \in \D\!\;.
\end{equation*}
If $\mathcal{ C}^*$ is finite, then $F$ is a finite Blaschke product.

Furthermore, $f=T \circ F$ for some automorphism $T$ of $\D$
if and only if
\begin{equation*}
\lim_{z \to w} \frac{|f'(z)|}{1-|f(z)|^2} \,
\frac{1-|F(z)|^2}{|F'(z)|}=1
\end{equation*}
for some $w \in \D$.
\end{theorem}

If  $\mathcal{ C}^*=\mathcal{ C} \not= \emptyset$, this gives  the  sharpening
\begin{equation*}
\frac{|f'(z)|}{1-|f(z)|^2} \le \frac{|F'(z)|}{1-|F(z)|^2} < \frac{1}{1-|z|^2}
\end{equation*}
for all $f \in \mathcal{F}_{\mathcal{C}}$
of the Schwarz--Pick inequality (\ref{eq:sp}), which is best possible in some sense.

\subsection{Related extremal problems in Hardy and Bergman spaces}

Let  $\mathcal{ C}$ be a sequence in $\D$, assume that $\mathcal{C}$ is the critical set of a  bounded analytic
function $f : \D \to \D$ and let  $N$ denote the multiplicity of the point $0$
in $\mathcal{C}$. Then according to Theorem \ref{thm:2}
 the maximal Blaschke product $F$ for $\mathcal{ C}$
normalized by $F(0)=0$ and $F^{(N+1)}(0)>0$ is the unique solution to the
extremal problem

\begin{center}
\fbox{\begin{minipage}{11cm}
{\begin{equation*}
 \max \big\{\Re f^{(N+1)}(0): f \in H^{\infty}\,, \, 
||f||_{\infty} \le 1 \text{ and } f'(z)=0 \text{ for } z \in\mathcal{ C} \big \}
\, .
\end{equation*} \vspace*{-3mm}}
\end{minipage}}
\end{center}

This extremal property of a maximal Blaschke product  is reminiscent 
of the well--known extremal property of 
\begin{itemize}
\item[(i)]
Blaschke products in the Hardy spaces $H^{\infty}$ and
$$H^p:=\left\{ f : \D \to \C \text{ analytic} \,: {||f||}_p<+ \infty \right\} \, , $$
where $1 \le p  < +\infty$ and
$$  {||f||}_p:=\left( \lim
    \limits_{r \to 1} 
\frac{1}{2\pi} \, \int
  \limits_{0}^{2\pi}  |f(r e^{it})|^p\,
dt  \right)^{1/p} \, ;$$
\end{itemize}
and 
\begin{itemize}
\item[(ii)] canonical divisors in the
(weighted) Bergman spaces  
$$ \mathcal{ A}_{\alpha}^p= \left\{
                 f:\D \to \C \text{ analytic} \,: {||f||}_{p,\alpha}<+ \infty \right\} \, , $$
where $-1<\alpha<+\infty$ and $1 \le  p<+\infty$ and
$$ {||f||}_{p,\alpha}:=\left(\frac{1}{\pi} \, \iint \limits_{\D} (1-|z|^2)^{\alpha} \, |f(z)|^p\,
d\sigma_z\right)^{1/p} \, .$$
\end{itemize} 
Note that in (i) and (ii) the prescribed data are not the critical points, but
the zeros.

\medskip
More precisely, let the sequence $ \mathcal{ C}=(z_j)$ in $\D$ be the zero set of
an $H^p$ function and let $N$ be  the multiplicity
of the point $0$ in $\mathcal{ C}$.  Then the (unique) solution to the extremal problem


\begin{center}
\fbox{\begin{minipage}{11cm}{\begin{equation*}
\max\big\{\Re f^{(N)}(0): f \in H^p,\,  ||f||_p \le 1 \text{ and } f(z)=0
\text{ for } z \in\mathcal{ C} \big\}
\end{equation*} \vspace*{-3mm}}
\end{minipage}}
\end{center}
is a Blaschke product $B$  with zero set $\mathcal{ C}$ which is normalized by
$B^{(N)}(0)>0$, see \cite[\S 5.1]{DS}.

In searching for an analogue of  Blaschke products for Bergman spaces, 
 Hedenmalm \cite{Hed1991} (see also \cite{DKSS1993,DKSS1994}) had the idea of posing an
 appropriate counterpart of the latter extremal problem for 
Bergman spaces. 
As before, let $\mathcal{ C}=(z_j)$
be a sequence in $\D$ where the point $0$ occurs $N$ times and assume that
$\mathcal{C}$ is the zero set of a function in $\mathcal{A}^p_{\alpha}$.
Then the extremal problem 


\begin{center}
\fbox{\begin{minipage}{11cm}{ \begin{equation*} 
 \max \big\{\Re f^{(N)}(0): f \in \mathcal{
  A}_{\alpha}^p\,, \,  || f||_{p,\alpha} \le 1 \text{ and } f(z)=0 \text{ for } z \in\mathcal{ C}
\big \}
\end{equation*}\vspace*{-3mm}}
\end{minipage}}
\end{center}
has a unique  extremal function $\mathcal{ G} \in \mathcal{A}^p_{\alpha}$,  which vanishes
precisely on $\mathcal{ C}$ and is normalized by $\mathcal{ G}^{(N)}(0)>0$. 
The function $\mathcal{G}$ is called the canonical divisor for $\mathcal{ C}$.
It plays a prominent r$\hat{\text{o}}$le in the modern theory of Bergman spaces.

\medskip

In summary, we have the following situation:

\smallskip
\renewcommand{\arraystretch}{1.5}
\begin{center}
\begin{tabular}{|l|c|l|}
\hline
{\bf prescribed data} & {\bf function space} & {\bf extremal function} \\[0mm]\hline
 critical set $\mathcal{ C}$  & $H^{\infty}$ & maximal Blaschke product for $\mathcal{ C}$\\[0mm]\hline
 zero set $\mathcal{ C}$       & $H^p$ & Blaschke product \\[0mm]\hline
zero set  $\mathcal{ C}$       & $\mathcal{ A}^p_{\alpha}$ & canonical divisor for
$\mathcal{ C}$ \\[0mm]
\hline
\end{tabular}
\end{center}

\medskip

In light of this strong analogy, 
one expects that maximal Blaschke products enjoy similar properties as finite
Blaschke products and canonical divisors. An example is their analytic continuability. It is
a familiar result that a Blaschke product 
has a holomorphic extension across every open arc of $\partial \D$ that does not
contain any limit point of its zero set, 
see \cite[Chapter II, Theorem 6.1]{Gar2007}.~The same is true for a canonical
divisor in the Bergman spaces $\mathcal{A}^p_0$. This
 was proved by
 Sundberg \cite{Sun1997} in 1997, who  improved earlier work
of Duren, Khavinson, Shapiro and Sundberg \cite{DKSS1993,DKSS1994} and Duren,
Khavinson and Shapiro \cite{DKS1996}.
Now following the model that critical points of maximal functions correspond
to the zeros of  Blaschke products and  canonical divisors respectively, one
hopes that a maximal Blaschke product has an analytic continuation across every open arc of
$\partial \D$ which does not meet any limit point of its critical set. This in fact
turns out to be true:

\begin{theorem}[Analytic continuability, \cite{KR2011}]\label{thm:analytic_continuation}
Let $F: \D \to \D$ be a maximal Blaschke product  with critical set $\mathcal{ C}$. Then $F$ has an
analytic continuation across each arc of $\partial \D$ which is free of limit
points of $\mathcal{ C}$.
In particular, the limit points of the critical set of $F$ coincide with the
limit points of the zero set of $F$.
\end{theorem}

Another rather strong property of finite Blaschke products is their semigroup
property with respect to composition. In contrast, the composition of two
 infinite Blaschke products does not need to be 
a Blaschke product (just consider destructible Blaschke products). However, in
the case of
maximal Blaschke products the following result holds.

\begin{theorem}[Semigroup property, \cite{KR2011}]
The set of maximal Blaschke products is closed under composition.
\end{theorem}

It would be interesting to get some information about the critical \textit{values} of
maximal Blaschke products and to explore the possibility of factorizing
maximal Blaschke products in a way similar to the recent extension of
Ritt's theorem for finite Blaschke products due to Ng and Wang (see \cite{NgWang2011}).

\subsection{Boundary behaviour of maximal Blaschke products}

We now shift attention to the boundary behaviour
of maximal Blaschke products. Ideally, one should be able to determine whether
a bounded analytic function $F : \D \to \D$
is a maximal Blaschke produkt either from the behaviour of
$$ \frac{|F'(z)|}{1-|F(z)|^2} \quad \text{ as } |z| \to 1$$
or from the behaviour of 
$$\int \limits_0^{2\pi} \log  \frac{|F'(r e^{it} )|}{1-|F(r e^{it})|^2} \, dt
\quad \text{ as } r \to 1 \, .$$

We only have some partial results in this connection and
we begin our account with the case of finite Blaschke products.

\begin{theorem}[Boundary behaviour of finite Blaschke products] \label{thm:boundary}
Let $I \subset \partial \D$ be some open arc and let
$f: \D \to \D$ be an analytic function.
Then the following statements are equivalent.
\begin{itemize}
\item[(a)] 
$ \lim \limits_{z \to \zeta} \left( 1-|z|^2 \right)
 \displaystyle \frac{|f'(z)|}{1-|f(z)|^2}=1 \qquad \text{ for every } 
\zeta \in I \, , $
\item[(b)]$ \lim \limits_{z \to \zeta} 
 \displaystyle \frac{|f'(z)|}{1-|f(z)|^2}=+\infty \qquad \text{ for every } 
\zeta \in I \, , $
\item[(c)]
$f$ has a holomorphic extension across the arc $I$ with $f(I) \subset \partial \D$.
\end{itemize}
In particular, if $I=\partial \D$, then $f$ is in either case a finite Blaschke product.
\end{theorem}

The equivalence of  conditions (a) and  (b) in Theorem \ref{thm:boundary}
for the special case $I=\partial \D$  is due to Heins \cite{Hei86}; the
general case is  proved in \cite{KRR06}. 
We now extend Theorem \ref{thm:boundary} beyond the class of finite Blaschke
products and start with the following  auxiliary result.

\begin{proposition}[see \cite{Kra2011a}] \label{lem:lemma2}
Let $f : \D \to \D$ be an analytic function and $I$ some subset of $\partial\D$.
\begin{itemize}
\item[(1)]
If 
\begin{equation*}
\angle \lim \limits_{z \to \zeta} \left( 1-|z|^2 \right)
 \displaystyle \frac{|f'(z)|}{1-|f(z)|^2}=1 \qquad \text{ for every } 
\zeta \in I \, , 
\end{equation*}
then  $f$ has a finite angular derivative\footnote{see \cite[p.~57]{Sha1993}.} at a.\!\;e.~$\zeta
  \in I$. In particular,  
\begin{equation*}
 \angle \lim \limits_{z \to \zeta} |f(z)|=1 \qquad \text{ for a.\!\;e. } \zeta \in  I
\, .
\end{equation*}

\item[(2)]
If $f$ has a finite angular derivative (and $\angle \lim_{z \to \zeta} |f(z)|=1$) at some $\zeta
  \in I$, then
\begin{equation*}
\angle \lim \limits_{z \to \zeta} \left( 1-|z|^2 \right)
 \displaystyle \frac{|f'(z)|}{1-|f(z)|^2}=1 \, .
\end{equation*}
\end{itemize}
\end{proposition}

In particular, when $I=\partial \D$, we obtain the following corollary.

\begin{corollary}[see \cite{Kra2011a}] \label{cor:5}
Let $f : \D \to \D$ be an analytic function. Then the following statements are
equivalent.
\begin{itemize}
\item[(a)] $\angle \lim \limits_{z \to \zeta} \left( 1-|z|^2 \right)
 \displaystyle \frac{|f'(z)|}{1-|f(z)|^2}=1$  for a.\!\;e.~$\zeta \in \partial
 \D$.
\item[(b)] $f$ is an inner function with finite angular derivative at almost
  every point of $\partial \D$.
\end{itemize}
\end{corollary}

We further note that
conditions (1) and (2) in Proposition \ref{lem:lemma2} do not complement each
other. Therefore we may ask if an analytic self--map $f$ of $\D$ which 
satisfies
\begin{equation*}
\angle\lim \limits_{z \to 1} \frac{|f'(z)|}{1-|f(z)|^2}\, (1-|z|^2)=1  
\end{equation*}
does have
an angular limit or even a finite angular derivative  at
$z=1$; this might then be viewed as a  converse of the
Julia--Wolff--Carath\'eodory theorem, see \cite[p.~57]{Sha1993}.

For maximal Blaschke products whose critical sets satisfy the Blaschke
condition one can show that condition (a) in Corollary \ref{cor:5} holds:

\begin{theorem}[see \cite{Kra2010}] \label{prop:1}
Let $\mathcal{ C}=(z_j)$ be a Blaschke sequence  in $\D$. 
\begin{itemize}
\item[(a)]
The maximal conformal pseudometric $\lambda_{max}(z) \, |dz|$
on $\D$ with  constant curvature $-4$ and zero set
$\mathcal{ C}$ satisfies 
\begin{equation*}
 \angle \lim \limits_{z \to \zeta} \ \frac{ \lambda_{max}(z)}{\lambda_{\D}(z)}=1
\quad \text{ for a.\!\;e. } \zeta \in \partial \D \, .
\end{equation*}
\item[(b)] \label{thm:6}
 Every maximal function
for $\mathcal{ C}$ has a finite angular derivative at almost
every point of $\partial \D$. 
\end{itemize}
\end{theorem}

\begin{remark} \label{satz:heinsproof}
The boundary behaviour of a maximal Blaschke product contains
useful information. For instance, it leads to a  quick proof of Theorem
\ref{satz:heins}. To see this, let $F: \D \to \D$
 be a maximal function for a finite sequence $\mathcal{ C}$ and $\lambda_{max}(z)\,
|dz|= (F^*\lambda_{\D})(z)\, |dz|$ be the maximal conformal metric with constant
curvature $-4$ and zero set $\mathcal{ C}$. By
Theorem \ref{thm:0}, the maximal function $F$ is an indestructible Blaschke product. 
Now, let $B$ be a finite Blaschke product with \textit{zero set} $\mathcal{ C}$.
Then 
\begin{equation*}
 \lambda(z)\, |dz| :=|B(z)| \, \lambda_{\D}(z)  \, |dz|
\end{equation*}
is a regular conformal pseudometric on $\D$ with curvature $-4/ |B(z)|^2\le -4$ and zero set
$\mathcal{ C}$.
Thus, by the maximality of $\lambda_{max}(z) \, |dz|$,
\begin{equation*}
\lambda(z) \le \lambda_{max}(z) \quad \text{ for } z \in \D 
\end{equation*}
and consequently
\begin{equation}\label{eq:ab1}
 |B(z)| \le \frac{\lambda_{max}(z)}{\lambda_{\D}(z)}  \quad \text{ for all } z \in \D  \, . 
\end{equation}
Since $B$ is a finite Blaschke product,  we deduce from (\ref{eq:ab1}) and the
Schwarz--Pick lemma (\ref{eq:sp}) that
\begin{equation*}
\lim_{z \to \zeta}
\frac{\lambda_{max}(z)}{\lambda_{\D}(z)}=\lim_{z \to \zeta} \frac{|F'(z)|}{1-|F(z)|^2}\, (1-|z|^2) =1 \quad \text{ for
  all } \zeta \in \partial \D \, .
\end{equation*}
Applying Theorem \ref{thm:boundary}  shows that $F$ is a finite
Blaschke product.
The branching order of $F$ is clearly $2\!\;n$. Thus,
according to the Riemann--Hurwitz formula, see \cite[p.~140]{For1999}, the
Blaschke product $F$ has degree $m=n+1$.

On the other hand, assume that $F$ is a finite Blaschke product of degree
$n+1$. Then, by Theorem \ref{thm:boundary},
$$ \lim_{z \to \zeta} \frac{|F'(z)|}{1-|F(z)|^2}\, (1-|z|^2) =1 \quad \text{ for
  all } \zeta \in \partial \D \, .$$
Theorem \ref{thm:3_11} below shows that $F$ is a maximal Blaschke product.
\end{remark}

The next result gives a sufficient condition for maximality
of a  Blaschke product $F$  in terms
of the boundary behaviour of the \textit{integral means} of the quantity
$$  (1-|z|^2) \frac{|F'(z)|}{1-|F(z)|^2}\,
=\frac{\lambda(z)}{\lambda_{\D}(z)} \, . $$

\begin{theorem}[see \cite{Kra2010}] \label{thm:3_11}
Let $\lambda(z)\, |dz|$ be a conformal pseudometric on $\D$ with constant curvature $-4$ and zero set $\mathcal{
  C}$ such that
\begin{equation}\label{eq:boundary2}
\lim_{ r\to 1} \,  \int \limits_0^{2\pi} \log
\frac{\lambda(re^{it})}{\lambda_{\D}(re^{it})}\,dt=0\, .
\end{equation}
Then $\lambda(z)\, |dz|$ is the maximal conformal pseudometric $\lambda_{max}(z)\, |dz|$ on $\D$
with constant curvature $-4$ and zero set $\mathcal{ C}$.
\end{theorem}

If $\mathcal{ C}$ is a  Blaschke sequence, then the corresponding maximal conformal pseudometric $\lambda_{max}(z)\, |dz|$ on $\D$
with constant curvature $-4$ and zero set $\mathcal{ C}$ satisfies
condition (\ref{eq:boundary2}):

\begin{theorem}[see \cite{Kra2011a}]
Let $\mathcal{ C}$ be a  Blaschke sequence  in $\D$. A conformal pseudometric
$\lambda(z)\, |dz|$ on $\D$ with constant curvature $-4$ and zero set $\mathcal{
  C}$ is the maximal conformal pseudometric $\lambda_{max}(z)\, |dz|$ on $\D$
with constant curvature $-4$ and zero set $\mathcal{ C}$ if and only if
(\ref{eq:boundary2}) holds.
\end{theorem}

We don't know whether this result is true for any sequence $\mathcal{ C}$ for
which there is a non--constant bounded analytic function with critical set
$\mathcal{ C}$.

\subsection{Heins' results on maximal functions}

For completeness,  we close this section with a discussion of Heins' results on maximal functions,
cf.~\cite[\S 25 \& \S 26]{Hei62}. A first observation is that every maximal
function is surjective. In fact more is true; a maximal function is ``locally''
surjective. Here is the  precise definition. 

\begin{definition}
Let $f: \D \to \D$ be an analytic function. A point $q  \in \D$ is called locally
omitted by $f$ provided that either $q \in \D \backslash f(\D)$ or else $ q
\in f(\D)$ and there exists a domain $\Omega$, $q \in \Omega$, such that for
some component $U$ of $f^{-1}(\Omega)$ the restriction of $f$ to $U$ omits $q$,
i.\!\;e.~$q \not\in f(U)$. 
\end{definition}

\begin{satz}[Heins \cite{Hei62}]\label{thm:omitted_value}
A maximal function has no locally omitted point.
\end{satz}

So far, the results about maximal functions leave an important question
unanswered, namely, how
to tell whether a given function in $H^{\infty}$ is a maximal function? Heins'
second result gives a topological sufficient criterion and 
provides therefore a source of ``examples'' of maximal functions. It is based
on the following concept.

\begin{definition}
A holomorphic function $f: \D \to \D$ is called locally of island type if $f$ is
onto and if for each $w \in \D$ there is an open disk $K(w)$ about
$w$ such that each component of $f^{-1}(K(w))$ is compactly contained in $\D$.
\end{definition}

Obviously, every surjective analytic self--map of $\D$ with constant finite
valence, that is, every finite Blaschke product is
locally of island type.

\begin{satz}[Heins \cite{Hei62}]\label{thm:heins_island}
Every function locally of island type is a maximal function. 
\end{satz}


\section{The Gauss curvature PDE and the Berger--Nirenberg problem}
\label{sec:gauss}

We return to a discussion of Theorem \ref{thm:main}. The results of Section
\ref{sec:2} (Theorem \ref{thm:perron} and Theorem  \ref{thm:0})
provide a proof of implication ``(a) $\Longrightarrow$ (b)'' in Theorem \ref{thm:main}.
In this section, we discuss implication
``(c) $\Longrightarrow$ (a)''.
The key idea is the following.

\begin{theorem}\label{thm:3}
Let $h : \D   \to \C$ be a non--constant holomorphic function with zero set
$\mathcal{C}$. Then the following statements are
equivalent.
\begin{itemize}
\item[(a)] There exists a holomorphic function $f : \D \to \D$ with critical
  set $\mathcal{C}$.
\item[(b)] There exists a $C^2$--solution $u: \D \to \R$ to the Gauss curvature equation
\begin{equation}\label{eq:gauss}
\Delta u = 4 \,|h(z)|^2 e^{2u}\, .
\end{equation}
\end{itemize}
\end{theorem}

Let us sketch a proof here. If $f : \D \to \D$ is a holomorphic function with critical
  set $\mathcal{C}$, then a quick computation shows that
$$ u(z):=\log \left( \frac{1}{|h(z)|} \frac{|f'(z)|}{1-|f(z)|^2} \right)$$
is a $C^2$--solution to (\ref{eq:gauss}). This proves ``(a) $\Longrightarrow$ (b)''.
Conversely, if there is  a $C^2$--solution $u: \D \to \R$ to the curvature equation (\ref{eq:gauss}),
then  $$\lambda(z)\, |dz |:=e^{u(z)}\, |dz|$$
 is a regular conformal metric with
curvature $-4 |h(z)|^2$ on $\D$. Hence, by  Theorem \ref{rem:liouville}, 
$$ u(z)=\log \left( \frac{1}{|h(z)|} \frac{|f'(z)|}{1-|f(z)|^2} \right)$$
with some analytic self--map $f$ of $\D$.
Thus the zero set of $h$ agrees with the critical set
of $f$.

In view of Theorem \ref{thm:3}, the task is now   to characterize those holomorphic functions $h: \D
\to \C$ for which  the PDE (\ref{eq:gauss}) has a solution.
In fact this problem is a special case of the 
Berger--Nirenberg problem from differential geometry: 

\medskip

{\bf Berger--Nirenberg problem:}\\
{\sl
Given a function $\kappa: R \to \R$ on a
Riemann surface $R$. Is there a conformal metric on $R$ with Gauss curvature
$\kappa$?} 

\medskip

The Berger--Nirenberg problem is well--understood for the projective plane,
see \cite{Mos1973} and has been extensively studied for compact Riemannian surfaces,
see \cite{Aub1998,Chang2004,Kaz1985,Str2005} as well as for the complex plane
\cite{Avi1986, CN1991, Ni1989}\footnote{These are just
some of the many references.}.  
However much less is known for proper domains $D$ of the complex plane,
see \cite{BK1986, HT92, KY1993}. In
this situation the Berger--Nirenberg problem
reduces  to the question if for a given function $k:D \to \R$ the Gauss curvature equation
\begin{equation}\label{eq:curvature1}
\Delta u= k(z)\, e^{2u}
\end{equation}
has a solution on $D$. We just note that $k$ is the
negative of the curvature $\kappa$ of the conformal metric $e^{u(z)}\, |dz|$.

In the next theorem  we give some necessary conditions as well as sufficient conditions
for the solvability of the Gauss curvature equation  (\ref{eq:curvature1})
only in terms of the curvature function $k$ and the domain $D$.

\begin{theorem}[see \cite{Kra2011a}]\label{thm:sol1}
Let $D$ be a bounded and regular domain\footnote{i.\!\;e.~there exists 
  Green's function $g_D$ for $D$ which vanishes
  continuously on $\partial D$.} and  let $k$  be a nonnegative locally H\"older
continuous function on $D$. 
\begin{itemize}
\item[(1)]
If for some (and therefore for every) $z_0 \in D$
\begin{equation*}
\iint \limits_{D} g_D(z_0, \xi)\, k(\xi) \, d\sigma_{\xi} < + \infty\, ,
\end{equation*}
then (\ref{eq:curvature1}) has a $C^2$--solution $u : D \to \R$, which is bounded from above.

\item[(2)]
If (\ref{eq:curvature1}) has a $C^2$--solution $u: D \to \R$ which is bounded from
below and if this solution has a harmonic majorant on $D$, then  
\begin{equation*}
\iint \limits_{D} g_D(z, \xi) \, k(\xi) \, d\sigma_{\xi} < + \infty
\end{equation*}
for all $z \in D$.
\item[(3)]
There exists a bounded $C^2$--solution 
$u:D \to \R$  to (\ref{eq:curvature1})  if and only if 
\begin{equation*}
\sup \limits_{z \in D} \iint \limits_{D} g_D(z, \xi) \,   k(\xi)\,   d\sigma_{\xi} < + \infty\, .
\end{equation*}
\end{itemize}
\end{theorem}

If we choose $D= \D$ and $z_0=0$, then $g_{\D}(0,\xi)=-\log|\xi|$. Hence as a consequence of
the inequality
\begin{equation*} 
 \frac{1-|\xi|^2}{2} \le \log \frac{1}{|\xi|} \le
\frac{1-|\xi|^2}{|\xi|}\, , \quad 0< |\xi| <1,
\end{equation*}
we obtain the following equivalent formulation of Theorem \ref{thm:sol1}.

\begin{corollary}[see \cite{Kra2011a}]\label{cor:solutions1}
Let $k$ be a nonnegative locally H\"older continuous function on $\D$. 
\begin{itemize}
\item[(1)]
If 
\begin{equation} \label{eq:sufcon}
\iint \limits_{\D} (1-|\xi|^2 )\, k(\xi) \, d\sigma_{\xi} < + \infty\, ,
\end{equation}
then (\ref{eq:curvature1}) has a $C^2$--solution $u : \D \to \R$, which is bounded from above.

\item[(2)]
If (\ref{eq:curvature1}) has a $C^2$--solution $u: \D \to \R$ which is bounded from
below and if this solution has a harmonic majorant on $\D$, then 
\begin{equation*}
\iint \limits_{\D} (1-|\xi|^2 )\, k(\xi) \, d\sigma_{\xi} < + \infty\, .
\end{equation*}
\item[(3)]
There exists a bounded $C^2$--solution $u: \D \to \R$  to (\ref{eq:curvature1}) if and only if 
\begin{equation*}
\sup \limits_{z \in \D} \iint \limits_{\D} \log\left| \frac{1- \overline{\xi} z}{z -
      \xi} \right|\,  k(\xi)\,   d\sigma_{\xi} < + \infty\, .
\end{equation*}
\end{itemize}
\end{corollary}
\smallskip

Both, Theorem \ref{thm:sol1} and  Corollary \ref{cor:solutions1}, are not best
possible, because (\ref{eq:curvature1}) may indeed have solutions, even if 
\begin{equation*}
\iint \limits_{D} g_{D}(z, \xi)\, k(\xi) \, d\sigma_{\xi} =+ \infty
\end{equation*}  
for some (and therefore for all) $z \in D$.
Here is an explicit example.

\begin{example} \label{ex:0}
For $\alpha \ge 3/2$ define
\begin{equation*}
h(z)= \frac{1}{(z-1)^{\alpha}}
\end{equation*}
for $z \in \D$ and set $k(z)=4\, |h(z)|^2$ for $z \in \D$.
Then an easy computation yields
\begin{equation*}
\iint \limits_{\D} (1-|z|^2 )\, k(z) \, d\sigma_{z}= +\infty
\end{equation*}
and a straightforward check shows that
 the function
\begin{equation*}
u_f(z):=  \log \left( \frac{1}{|h(z)|} \, \, \frac{|f'(z)|}{1-|f(z)|^2}
\right) 
\end{equation*}
is a solution to
(\ref{eq:curvature1}) on $\D$ for every locally univalent analytic function $f: \D \to \D$.
\end{example}

\begin{remark}
The sufficient condition (\ref{eq:sufcon})
improves earlier results of Kalka $\&$ Yang in \cite{KY1993}. In fact, Kalka $\&$
Yang give  explicit examples for the function  $k$ which tend to $+\infty$ at the
boundary of $\D$ such that  the existence of a solution to
(\ref{eq:curvature1}) can  be
guaranteed. All these examples are radially symmetric  and satisfy
(\ref{eq:sufcon}). 
Kalka $\&$ Yang also supplement their existence results by nonexistence results.  
They find explicit lower bounds for the function $k$ in terms of 
 radially symmetric functions which grow to $+\infty$ at the boundary of $\D$,
such that
(\ref{eq:curvature1}) has no solution.

We wish to emphasize that the necessary conditions and the sufficient conditions
for the solvability of the curvature equation (\ref{eq:curvature1}) of Kalka \& Yang
do not complement each other. In particular, the case when 
the function $k$  oscillates is not covered.
For the proof of their nonexistence results 
Kalka $\&$ Yang needed to use Yau's celebrated Maximum Principle
\cite{Yau1975, Yau1978}, which is an extremely
powerful tool. In \cite{Kra2010},  an almost elementary proof of these
nonexistence results is given, which has the additional
advantage that  Ahlfors' type lemmas for  conformal metrics with
variable curvature and explicit formulas for the corresponding maximal
conformal metrics are obtained. In \cite{Kra2010} the nonexistence
theorems of Kalka \& Yang are further extended by allowing the function $k$ to oscillate. 
\end{remark}

It turns out that
although condition (\ref{eq:sufcon}) is not necessary for the existence of a
solution to (\ref{eq:curvature1}) it is strong enough to deduce a necessary and
sufficient condition for the
solvability of the Gauss  curvature equation of the
particular form (\ref{eq:gauss}):

\begin{theorem}[see \cite{Kra2011a}] \label{thm:k_holo}
Let $h:\D \to \C$ be a holomorphic function. Then the Gauss curvature equation
(\ref{eq:gauss}) has a solution if and only if $h$ has
a representation as a product of an $\mathcal{ A}_1^2$ function and a nonvanishing
analytic function.
\end{theorem}

Note that  Theorem \ref{thm:3} combined with Theorem \ref{thm:k_holo}  shows that the class of all holomorphic functions
$h: \D \to \C$ whose zero sets  coincide with the critical sets of the class
of bounded analytic function is exactly the Bergman space $\mathcal{A}^2_1$.
This proves implication ``(c) $\Longrightarrow$ (a)'' in Theorem \ref{thm:main}.

A further remark is that
Theorem \ref{thm:sol1} (c)  characterizes those curvature functions $k$ 
for which (\ref{eq:curvature1}) has at least one
bounded solution. For the case of the unit disk $\D$, this result can be
stated as follows.

\begin{theorem}\label{rem:boundedsol}
Let $\varphi: \D \to \C$  be analytic and $k(z)=4 \, |\varphi'(z)|^2$.
Then there exists a bounded solution to the Gauss curvature equation (\ref{eq:curvature1}) if and only if
$\varphi \in \BMOA$, where 
\begin{equation*}
\BMOA=\left\{ \varphi: \D \to \C \text{ analytic} \, : \,  \sup_{ z \in \D} \iint \limits_{\D} g_{\D}(z,
  \xi)\, |\varphi'(\xi)|^2 \, d\sigma_{\xi} < + \infty \right \}
\end{equation*}
is the space of analytic functions of bounded
mean oscillation on $\D$, see \cite[p.~314]{BF2008}.
\end{theorem}

Finally, we  note that in Theorem \ref{thm:sol1} and Corollary
\ref{cor:solutions1}, condition (1) does not imply condition (3). 
The Gauss curvature equation (\ref{eq:curvature1})   may indeed
have solutions none of which is bounded. For example, choose
$\varphi \in H^2\backslash BMOA$ and set $k(z)=4\, |\varphi'(z)|^2$. Then, according to Theorem \ref{rem:boundedsol}, every solution to
(\ref{eq:curvature1}) must be unbounded.

The following result of Heins adds another
item to the list of equivalent statements in Theorem \ref{thm:main}.

\begin{satz}[Heins \cite{Hei62}] \label{satz:heins2}
Let $\mathcal{C}=(z_j)$ be a sequence in $\D$. Then the following conditions are equivalent.
\begin{itemize}
\item[(a)]  There is an analytic function $f : \D \to \D$  with critical set $(z_j)$.
\item[(b)] There is a function in the Nevanlinna class $\mathcal{ N}$ with critical set $(z_j)$.
\end{itemize}
\end{satz}

Here,  a function $f$  analytic in $\D$ is said to belong to the
Nevanlinna class $\mathcal{ N}$ if the integrals $$\int \limits_0^{2\pi} \log^+
|f(re^{it})|\, dt$$ remain bounded as $r \to 1$.

Let us indicate how the results of the present survey allow a quick proof of
Theorem~\ref{satz:heins2}.

\begin{proof}
(b) $\Rightarrow$ (a):  Let $\varphi \in \mathcal{ N}$. Then $\varphi=
\varphi_1/\varphi_2$ is the quotient of two 
analytic self--maps of $\D$, see for instance
 \cite[Theorem 2.1]{Dur2000}. W.\!\;l.\!\;o.\!\;g.~we may assume $\varphi_2$ is zerofree. Differentiation of $\varphi$ yields
\begin{equation*}
\varphi'(z)=\frac{1}{\varphi_2(z)^2} \, \big(
 \varphi_1'(z) \, \varphi_2 (z)- \varphi_1(z)\, \varphi_2'(z)
\big)\, .
\end{equation*}
Since $\varphi'_1, \varphi_2' \in \mathcal{ A}_1^2$ and $\mathcal{ A}_1^2$ is a vector
space,   it follows that the function $\varphi_1' \, \varphi_2 - \varphi_1\, \varphi_2'$
belongs to $\mathcal{ A}_1^2$. Thus Theorem \ref{thm:k_holo} ensures the existence of a
solution $u : \D \to \R$ to 
\begin{equation*}
\Delta u= \left| \varphi'(z) \right|^2 \, e^{2u}\, .
\end{equation*}
Hence Theorem \ref{thm:3} gives the desired result.
\end{proof}

The results of  this section about the solvability of the Gauss
curvature equation (\ref{eq:curvature1}) do have further
consequences for the critical sets of bounded analytic functions.
For instance, one can give an answer  to a question of Heins \cite[\S
31]{Hei62}.
In order to state Heins' question properly, 
we recall the well--known Jensen formula which connects the rate of growth
of an analytic function with  the density of its zeros.
Thus the restriction on the growth of the derivative of an analytic self--map
of $\D$ imposed by the Schwarz--Pick lemma (\ref{eq:sp}) forces
 an upper bound
for the  number of critical points of non--constant analytic self--maps of
$\D$. More precisely,
\begin{equation}\label{eq:anzahl1}
 \sup_{ {f \in
     H^{\infty}}\atop{||f||_{\infty}\le 1, \, f \not \equiv \text{const.}}}  \left( \limsup_{r \to 1} \,   \frac{N(r;f')}{\log\frac{1}{1-r^2}} \right)\le 1\, ,
\end{equation}
where 
$$ N(r;f'):=\int \limits_0^1
   \frac{ n(t;f')}{t} \, dt$$ and
$n(r;f')$ denotes the number of zeros of $f'$ counted with multiplicity in the
disk $\D_r:=\{ z \in \C\, : |z|<r\}$, $0<r<1$. 

\medskip
Heins showed that equality holds in (\ref{eq:anzahl1}), cf.~\cite{Hei61}.
 More precisely, using his solution of the
Schwarz--Picard problem,  he proved that 
there exists for every $p=2,3, \ldots$ a non--constant analytic function $f_p: \D \to \D$
such that
 \begin{equation*}
\limsup_{r \to 1 } \frac{N(r;f_p')}{\log \left(\frac{1}{1-r^2}\right)} =
\frac{2p-3}{2p-2}\, .
\end{equation*}
Letting $p \to + \infty$ shows that (\ref{eq:anzahl1}) is best possible. 
In this way, Heins  \cite[\S
31]{Hei62} was led to ask whether
there is a bounded  analytic function which realizes the supremum in 
(\ref{eq:anzahl1}).
This question is answered in our next theorem -- even with some
additional information.

\begin{theorem} \label{thm:heinsi}
Let $\beta \in [0,1]$. Then there exists a non--constant analytic function (and
even a maximal Blaschke product) $f : \D \to \D$ such that 
\begin{equation*} \limsup \limits_{r \to 1} \frac{N(r;f')}{\log \left(\frac{1}{1-r^2}\right)} =
\beta \, .
 \end{equation*}
\end{theorem}

For the proof of Theorem \ref{thm:heinsi} we refer the reader to \cite{Kra2010}.

We close with the following remark.

\begin{remark} \label{rem:final}
With the help of the Riemann mapping theorem,
the results of this paper about critical sets of  bounded analytic functions $f
: \D \to \D$ can easily be transferred to the class $H^{\infty}(\Omega)$
of bounded analytic functions $f : \Omega \to \D$, when $\Omega{\not=}\C$ is
a simply connected domain. The critical sets of bounded analytic functions
on multiply connected domains are much more difficult to fathom.
\end{remark}

\section*{Acknowledgement}

This paper is based on lectures given at the workshop on {\it Blaschke
products and their Applications} (Fields Institute, Toronto, July 25--29, 2011).
The authors would like to thank the organizers of this workshop,
Javad Mashreghi and Emmanuel Fricain, 
as well as  the Fields Institute and its staff, for their
generous support and hospitality.

\end{document}